\documentstyle[12pt]{article}
\newcommand{\card}{\mathop{\rm card}\limits}
\newcommand{\const}{\mathop{\rm const}\limits}
\newcommand{\Law}{\mathop{\rm Law}\limits}
\newcommand{\mod}{\mathop{\rm mod}\limits}
\newcommand{\Var}{\mathop{\rm Var}\limits}

\newcommand{\grad}{\mathop{\rm grad}\limits}

\newcommand{\supp}{\mathop{\rm supp}\limits}
\newcommand{\cov}{\mathop{\rm cov}\limits}
\newcommand{\extr}{\mathop{\rm extr}\limits}
\newcommand{\edim}{\mathop{\rm edim}\limits}

\begin{document}

\begin{center}

{\bf MULTIDIMENSIONAL PROBABILISTIC REARRANGEMENT} \\

\vspace{3mm}

{\bf INVARIANT SPACES: A NEW APPROACH.}\\

\vspace{4mm}

{\bf E. Ostrovsky and L. Sirota}\\

\vspace{4mm}

{\it Department of Mathematics and Statistics, Bar-Ilan University,
59200, Ramat Gan.}\\
e \ - \ mails: eugostrovsky@list.ru; \ sirota3@bezeqint.net \\

\vspace{4mm}

                 {\it Abstract.}\par
\end{center}

 We introduce and investigate in this paper a new convenient method of introduction
of a norm  in the multidimensional rearrangement probability invariant space. \par

 \vspace{4mm}

           {\it Key words and phrases:}\par

  Martingale and martingale differences, exponential and moment
 estimations, rearrangement invariant, moment and Grand Lebesgue-Riesz spaces of random variables,
 multidimensional rearrangement invariant (m.d.r.i.) spaces, fundamental function,  tail of
 distribution, polar, extremal points, random processes and fields, slowly and regular varying
 functions, sub-gaussian and pre-gaussian random variables, vectors and processes (fields),
 entropy and entropy integral, Young-Fenchel or Legendre transform.\par

 \vspace{4mm}

\par

\vspace{4mm}

\section{Introduction. Notations.}

\vspace{3mm}

 Let $ (\Omega,F,{\bf P} ) $ be a probability space, $ \Omega = \{\omega\},
\ (X, ||\cdot||X) $  be any rearrangement (on the other terms, symmetrical)
Banach function space over $ (\Omega,F,{\bf P} ) $ in the terminology of the
classical book \cite{Bennet1}, chapters 1,2;  see also \cite{Krein1}; $ d=2,3,4, \ldots, $

$$
S(d) = \left\{b, b \in R^d, |b|_2 := \left[\sum_{j=1}^d b^2(j) \right]^{1/2} \le 1 \right\}
$$
be the Euclidean unit ball with the center at the origin, $ B  $ be any separated complete
closed subset of the set $ S(d). $ \par
 The last imply that for arbitrary $ x \in R^d $ there exist two vectors
$  b_1,b_2 \in B $ for which $ (b_1,x) \ne (b_2,x). $  \par

{\bf Definition 1.} A Banach function space $ (X^{(d)},B) $ with the norm $ || \cdot||(X^{(d)},B) $
consists by definition on all the $ d- $ dimensional
random vectors (r.v.) $ \xi = \vec{\xi} = \{\xi(1), \xi(2), \ldots, \xi(d) \} $  with finite norm

$$
||\vec{\xi}||(X^{(d)},B) = ||\xi||(X^{(d)},B) \stackrel{def}{=}
\sup_{b \in B} ||(\xi,b)||X \eqno(1.1)
$$
is called {\it multidimensional rearrangement invariant (m.d.r.i.) space based on the space } $ X $
{\it (Khintchin's version) and the set } $ B. $ \par
 Here and furthermore $(b_1,b_2) $  denotes as customary the inner (scalar) product of two $ d- $  dimensional
vectors $ b_1,b_2: $

$$
(b_1,b_2) = \sum_{j=1}^d b_1(j) b_2(j), \ b_1 = \{b_1(j)\}, \ b_2 = \{b_2(j)\}, \ j=1,2,\ldots,d.
$$
 Obviously, instead the whole set $ B $ in (1.1) it may be stand the set of extremal points one. \par

 We will write for brevity in the case $ B = S(d) $

$$
||\vec{\xi}||(X^{(d)},S(d)) = ||\xi||(X^{(d)},S(d)) = ||\xi||X^{(d)}. \eqno(1.2)
$$

\vspace{3mm}

 Evidently, the space $ (X^{(d)},B) $ is rearrangement invariant in the ordinary sense: the norm
$ ||\vec{\xi}||(X^{(d)},B) $ dependent only on the distribution of the random vector $ \xi =\vec{\xi}:$

$$
\mu_{\xi}(A) = {\bf P}(\xi \in A),
$$
where $ A $ is arbitrary Borel set in the whole space $ R^d. $ \par
\vspace{4mm}

 Another notations: the $ L(p) = L_p = L_p(\Omega) $ is the classical Lebesgue-Riesz space
consisting on all the random variables $ \{\eta\} $
defined on the source probability space with finite norm

$$
|\eta|_p \stackrel{def}{=} \left[{\bf E} |\eta|^p \right]^{1/p}, \ p \ge 1.
$$

 For instance, the r.i. space $ (X, ||\cdot||X) $ may be the classical $ L_p(\Omega) $ space, Lorentz,
Marzinkievicz, Orlicz and so one spaces. We recall now for readers convenience briefly the definitions
and simple properties of some new r.i. spaces, namely, Grand Lebesgue spaces and Banach spaces $ \Phi(\phi) $
of random variables with exponentially decreasing tails of distribution (pre-gaussian and sub-gaussian
spaces).\par

 For $a = \const > 1$ let $\psi =
\psi(p),$ $p \in [1,a) $ be a continuous positive
function such that there exists a limits (finite or not)
$ \psi(a - 0)$ with conditions $ \inf_{p \in (1,a)} > 0. $
 We will denote the set of all these functions
 as $ \Psi(a). $ \par
 The Grand Lebesgue Space (in notation GLS) $  G(\psi; a) =
 G(\psi) $ is the space of all measurable
functions (random variables) $ \ \eta: \Omega \to R \ $ endowed with the norm
$$
||\eta||G(\psi) \stackrel{def}{=}\sup_{p \in (a,b)}
\left[ \frac{ |\eta|_p}{\psi(p)} \right]
$$
if it is finite.\par

These spaces are investigated, e.g. in the articles  \cite{Buldygin1}, \cite{Fiorenza1},
 \cite{Fiorenza2}, \cite{Fiorenza3}, \cite{Fiorenza4}, \cite{Ivaniec1}, \cite{Ivaniec2},
 \cite{Kozatchenko1}, \cite{Liflyand1}, \cite{Ostrovsky2}, \cite{Ostrovsky3}, \cite{Ostrovsky4},
 \cite{Ostrovsky5} etc. Notice that in the articles
\cite{Fiorenza2}, \cite{Fiorenza3}, \cite{Fiorenza4},
\cite{Ivaniec1}, \cite{Ivaniec2}, \cite{Liflyand1}, \cite{Ostrovsky3}, \cite{Ostrovsky4}
it was considered more general case when the measure $ {\bf P} $ may be unbounded. \par

  For instance, in the case when  $ a <  \infty, \beta = \const > 0 $ and

 $$
 \psi(p) = \psi(a, \beta; p) = (p - a)^{-\beta},  1 \le p < a;
 $$
we will denote the correspondent $ G(\psi) $ space by  $ G(a, \beta);  $ it
is not trivial, non-reflexive, non-separable etc.  \par
  Note that by virtue of Lyapunov inequality in the case $ \beta = 0, $ i.e.  when the  function
$ \psi(a;p) $ is bounded inside the {\it closed} interval $ p \in [1,a] $ and is infinite in the
exterior one, the space $ G(a, \beta)$ is equivalent to the classical Lebesgue-Riesz space $ L(a). $ \par
 In the case when $ a = \infty $ we need to take $ \gamma = \const > 0 $ and define

$$
\psi(p) = \psi(\gamma; p) = p^{m}, p \ge 1, \ m = \const > 0.
$$

 We will denote for simplicity these spaces as $ G(m). $ \par
We obtain a more general case considering the $ \psi $  function of a view

$$
\psi(p) = p^m \ L(p), \ p \ge 1,
$$
where the function $ L = L(p) $ is positive continuous slowly varying as $ p \to \infty $
function, so that the function $ \psi(p) $ is regular varying as $ p \to \infty. $ \par
 The following {\it sub-examples} are used in practice, see \cite{Liflyand1}, \cite{Ostrovsky5}:

 $$
 \psi(p) = p^m \ \log^{\gamma}(p+1) \ L(\log p).
 $$

  Further, let $ \phi = \phi(\lambda), \lambda \in (-\lambda_0, \lambda_0), \ \lambda_0 =
 \const \in (0, \infty] $ be some even strong convex which takes positive values for
 positive arguments twice continuous differentiable function, such that
$$
 \phi(0) = 0, \ \phi^{//}(0) \in(0,\infty), \ \lim_{\lambda \to \lambda_0}
 \phi(\lambda)/\lambda = \infty.
$$
 We denote the set of all these function as $ \Phi; \ \Phi =
\{ \phi(\cdot) \}. $ \par
 We will say that the {\it centered} random variable (r.v) $ \eta = \eta(\omega) $
belongs to the space $ \Phi(\phi), $ if there exists some non-negative constant
$ \tau \ge 0 $ such that

$$
\forall \lambda \in (-\lambda_0, \lambda_0) \ \Rightarrow
{\bf E} \exp(\lambda \eta) \le \exp[ \phi(\lambda \ \tau) ].
$$
 The minimal value $ \tau $ satisfying this inequality  is called a $ \Phi(\phi) \ $ norm
of the variable $ \xi, $ write
 $$
 ||\eta||\Phi(\phi) = \inf \{ \tau, \ \tau > 0: \ \forall \lambda \ \Rightarrow
 {\bf E}\exp(\lambda \xi) \le \exp(\phi(\lambda \ \tau)) \}.
 $$
 This spaces are very convenient for the investigation of the r.v. having a
exponential decreasing tail of distribution, for instance, for investigation of the limit
theorem, the exponential bounds of distribution for sums of random variables,
non-asymptotical properties, problem of continuous of random fields,
study of Central Limit Theorem (CLT) in the Banach spaces etc.\par

  The space $ \Phi(\phi) $ with respect to the norm $ || \cdot ||\Phi(\phi) $ and
ordinary operations is a Banach space which is isomorphic to the subspace
consisted on all the centered variables of Orlicz's space $ (\Omega,F,{\bf P}), N(\cdot) $
with $ N \ - $ function
$$
N(u) = \exp(\phi^*(u)) - 1, \ \phi^*(u) = \sup_{\lambda} (\lambda u -
\phi(\lambda)),
$$
see  \cite{Kozatchenko1}, \cite{Ostrovsky2}, chapter 1.\par
 The transform $ \phi \to \phi^* $ is called Young-Fenchel transform. The proof of considered
assertion used the properties of saddle-point method and theorem of Fenchel-Moraux:
$$
\phi^{**} = \phi.
$$

 The next facts about the $ B(\phi) $ spaces are proved in  \cite{Buldygin1},
 \cite{Kozatchenko1}, \cite{Ostrovsky2}, chapter 1, \cite{Ostrovsky3}:

$$
{\bf 1.} \ \eta \in B(\phi) \Leftrightarrow {\bf E } \eta = 0, \ {\bf and} \ \exists C = \const > 0,
$$

$$
T(\eta,x) \le \exp(-\phi^*(Cx)), x \ge 0,
$$
where $ T(\eta,x)$ denotes in this article the {\it tail} of
distribution of the r.v. $ \eta: $

$$
T(\eta,x) = \max \left( {\bf P}(\eta > x), \ {\bf P}(\eta < - x) \right),
\ x \ge 0,
$$
and this estimation is in general case asymptotically exact. \par
 More exactly, if $ \lambda_0 = \infty, $ then the following implication holds:

$$
\lim_{\lambda \to \infty} \phi^{-1}(\log {\bf E} \exp(\lambda \eta))/\lambda =
K \in (0, \infty)
$$
if and only if

$$
\lim_{x \to \infty} (\phi^*)^{-1}( |\log T(\eta,x)| )/x = 1/K.
$$
 Here and further $ f^{-1}(\cdot) $ denotes the inverse function to the
function $ f $ on the left-side half-line $ (C, \infty). $ \par
 {\bf 2.} We define $ \psi(p) = p/\phi^{-1}(p), \ p \ge 2. $
 Let us introduce a new norm (the so-called "moment norm")
on the set of r.v. defined in our probability space by the following way: the space
$ G(\psi) $ consist by definition on all the centered r.v. with finite norm

$$
||\eta||G(\psi) \stackrel{def}{=} \sup_{p \ge 2} [ |\eta|_p/\psi(p)].
$$

 It is proved in particular that the spaces $ \Phi(\phi) $ and $ G(\psi) $ coincides:$ \Phi(\phi) =
G(\psi) $ (set equality) and both
the norm $ ||\cdot||\Phi(\phi) $ and $ ||\cdot|| $ are equivalent: $ \exists C_1 =
C_1(\phi), C_2 = C_2(\phi) = \const \in (0,\infty), \ \forall \eta \in \Phi(\phi) $

$$
||\eta||G(\psi) \le C_1 \ ||\eta||B(\phi) \le C_2 \ ||\xi||G(\psi).
$$

{\bf 3.} This definition is correct still for the non-centered random
variables $ \eta.$ If for some non-zero r.v. $ \eta \ $ we have $ ||\eta||G(\psi) < \infty, $
then for all positive values $ u $

$$
{\bf P}(|\eta| > u) \le 2 \ \exp \left( - u/(C_3 \ ||\eta||G(\psi)) \right).
$$
and conversely if a r.v. $ \xi $ satisfies the last inequality, then $ ||\xi||G(\psi) <
\infty. $ \par

 Let $ \eta: \Omega \to R $ be some r.v. such that

 $$
 \eta \in \cup_{p > 1} L(p).
 $$
We can then introduce the non-trivial function $ \psi_{\eta}(p) $ as follows:

$$
\psi_{\eta}(p) \stackrel{def}{=} |\eta|_p.
$$
This choosing of the function $ \psi_g(\cdot) $ will be called {\it natural
choosing. }\par
 Analogously, if the centered (zero mean)  r.v. $ \eta $ satisfies the Kramer's condition

$$
\exists \mu \in (0, \infty), \ T(\eta, \ x) \le \exp(-\mu \ x),
\ x \ge 0,
$$
the function $ \phi(\cdot) = \phi_{\eta}(\lambda) $ may be "constructive" introduced by the formula
$$
\phi(\lambda) = \phi_0(\lambda) \stackrel{def}{=} \log
 {\bf E} \exp(\lambda \eta),
$$
 if obviously the  centered r.v. $ \eta $ satisfies the Kramer's condition:
$$
\exists \mu \in (0, \infty), \ T(\eta, \ x) \le \exp(-\mu \ x),
\ x \ge 0.
$$
 We will call also in this case the function $ \phi(\lambda) = \phi_{\eta}(\lambda) $
a {\it natural } function for the r.v. $ \eta. $ \par

\vspace{4mm}
 The letters $ C, C_k, C_k(\cdot),  k=1,2,\ldots $ with or without subscript  will denote a
finite positive non essential constants, not necessarily  the same at each appearance. \par

\vspace{4mm}

{\bf Our aim is investigation of properties of m.d.r.i. spaces: some inequalities, conjugate and
associate spaces, multidimensional moment and tail inequalities, estimation of normed sums of
independent random vectors or martingale differences random vectors, finding of sufficient conditions
for continuity of vector random processes and fields etc. } \par

\vspace{4mm}

 {\bf Example 1.1.} Let $ \xi = \vec{\xi} = \{\xi(1), \xi(2), \ldots, \xi(d) \} $ be a centered random
vector  with finite covariation matrix

$$
R(\xi) = \{ R_{i,j}(\xi)\}, \ R_{i,j}(\xi) = \cov(\xi(i), \xi(j)) = {\bf E} \xi(i)\cdot \xi(j), \
i,j = 1,2,\ldots,d.
$$
 If we choose $ X = L_2= L_2(\Omega), $ then

 $$
 ||\xi||^2(L_2^{(d)}) = \max_{b \in S(d)} {\bf E}(\xi,b)^2 =  \max_{b \in S(d)}
 \sum \sum_{i,j = 1,2,\ldots,d} R_{i,j}(\xi) b(i) b(j) = \lambda_{\max}(R(\xi)),
 $$
where  $ \lambda_{\max}(R(\xi))$ denotes the maximal eigen value of the matrix $ \lambda_{\max}(R(\xi)).$
Therefore,
$$
||\xi||(L_2^{(d)}) =   [\lambda_{\max}(R(\xi)) ]^{1/2}.
$$

 {\bf Example 1.2.} Let $ \xi = \vec{\xi} = \{\xi(1), \xi(2), \ldots, \xi(d) \} $ be a centered
{\it Gaussian} random vector  with covariation matrix $ R(\xi) $ and $ (X, ||\cdot||X) $ be arbitrary r.i.
space. As long as the random variable $ (\xi,b) $ has the mean zero Gaussian distribution with variance
$ \Var((\xi,b)) = (R(\xi)b,b), $ we conclude

$$
||\xi||(X^{(d)}) = [\lambda_{\max}(R(\xi)) ]^{1/2} \cdot ||\tau||X,
$$
where $ \tau $ has a standard Gaussian distribution.\par

 If $ X = L_p(\Omega), $ then

$$
||\tau||X = ||\tau||L_p(\Omega) = \sqrt{2} \ \pi^{-1/(2p)} \ \Gamma^{1/p}((p+1)/2),
$$
where $ \Gamma(\cdot) $ is Gamma function.\par
{\bf Example 1.3.} We can choose instead the $ (X, ||\cdot||X) $ many other r.i. spaces over
our probability space $ (\Omega,F,{\bf P} ), $  for instance, the Marzinkievicz $ M $ spaces,
Lorentz spaces $ L, $  Orlicz's space $ N $ with correspondent $ N-$ function $ N = N(u) $ and so ones. \par
 The detail investigation  of these spaces see in the classical monographs \cite{Bennet1}, chapters 1,2;
 \cite{Krein1}, chapters 1,2; see also  \cite{Astashkin1}.\par

{\bf Example 1.4.} In the capacity  of the space $ X $ it may be represented the space
$  G(\nu;r), $ which consist, by definition, on all the  r.v. with finite norm

$$
||\xi||G(\nu;r) \stackrel{def}{=} \sup_{p \in (2,r)} [|\xi|_p/\nu(p)], \ |\xi|_p :=
{\bf E}^{1/p} |\xi|^p.
$$
  Here $ r = \const > 2, \ \nu(\cdot) $ is some continuous positive on the
 {\it semi-open} interval  $ [1,r) $ function such that

$$
\inf_{p \in (2,r)} \nu(p) > 0, \ \nu(p) = \infty, \ p > r.
$$

 We will denote

$$
 \supp (\nu) \stackrel{def}{=} = \{p: \nu(p) < \infty \}.
$$
 {\it Sub-example:}

$$
\nu(p) = (r-p)^{-\gamma} \ L(1/(r-p)), \ r = \const >1,  \ 1 \le p < r, \ \gamma = \const \ge 0.
$$
where as before $  L = L(u) $ is positive continuous slowly varying as $ u \to \infty $ function.
 About applications of these spaces see  \cite{Ostrovsky6}. \par
 Other examples see in the section 5. \par

 The paper is organized as follows. In the section 2 are investigated simple properties of
introduces spaces.  In the section 3 we study conjugate and associate spaces. In the next section
we obtain some multidimensional tail inequalities for random vectors belonging to these spaces. \par
 In the fifth section we formulate and prove a moment inequalities for sums of random vectors.
The next section contains some information about fundamental function for m.d.r.i. spaces.
The $ 7^{th} $ section is devoted to a particular but important case of these spaces, namely, the
so-called multidimensional spaces of random vectors of subgaussian  and pre-gaussian type.
Multidimensional random fields are considered in the $ 8^{th} $ section.\par
   The last section contains some concluding remarks.\par

\vspace{4mm}
 \section{Simple properties of multidimensional rearrangement invariant (m.d.r.i.) spaces.}
\vspace{4mm}

{\bf Proposition 2.1.} Let $ \xi = \vec{\xi} = \{\xi(1), \xi(2), \ldots, \xi(d) \} $ be a
r.v. from the space $ X^{(d)}, $ then

$$
\max_{i=1,2,\ldots,d}||\xi(i)||X \le  ||\vec{\xi}||(X^{(d)})  \le
\sum_{i=1}^d||\xi(i)||X; \eqno(2.1a)
$$

$$
C_1 \max_{i=1,2,\ldots,d}||\xi(i)||X \le  ||\vec{\xi}||(X^{(d)},B)  \le
C_2 \sum_{i=1}^d||\xi(i)||X. \eqno(2.1b)
$$

{\bf Proof.}  As long as $ b \in S(d), \ |b(i)| \le 1, $ therefore

$$
||\vec{\xi}||X^{(d)} = \max_{b \in B} |(\xi,b)|X = \max_{b \in B} | \sum_{i=1}^d b(i) \xi(i)|X \le
\sum_{i=1}^d |\xi(i)|X.
$$
 We prove the right hand side of inequality (2.1a). \par
  Further, we have choosing the vector $ b $ of a view $ b = (0,0, \ldots, 1, 0,\ldots,0) =: e(k), $
where "1" stands on the place $ k, \ k=1,2,\ldots,d: $

$$
||\vec{\xi}||X^{(d)} \ge ||\xi(k)||X,
$$
hence

$$
||\vec{\xi}||X^{(d)} \ge \max_{k=1,2,\ldots,d} ||\xi(k)||X.
$$
 The assertion (2.1b) may be proved analogously.\par
 As a little consequences:\par
\vspace{3mm}
{\bf Proposition 2.2.} \par
{\bf A.} If the space $ (X, ||\cdot||X) $ is separable, then $ (X^{(d)},B) $ is separable. \par
{\bf B.} If the space $ (X, ||\cdot||X) $ is reflexive, then $ (X^{(d)},B) $ is reflexive, as well as. \par
\vspace{3mm}
 Note that an inverse conclusion is obvious. \par
 \vspace{3mm}
{\bf Proposition 2.3.} \par
 It is easy to formulate the criterions for convergence of the sequences in these spaces and compactness
of the sets.  Indeed,   convergence of the sequences is equivalent to coordinate-wise convergence;
for the compactness of the sets are true the classical features of Kolmogorov and Riesz. \par

\vspace{3mm}

\section{Conjugate and associate spaces.}

\vspace{3mm}

 We denote the {\it conjugate, or dual} space to the space $ (X, ||\cdot||X) $  as $ (X^*, ||\cdot||X^*) $
and {\it associate } space to the space $ (X, ||\cdot||X) $  as $ (X', ||\cdot||X'). $ \par

 By definition, arbitrary linear continuous functional on the space $ X^{(d)}, $ in the other words,
element of associate space of a view

 $$
 l_g(\xi) = \sum_{i=1}^d \int_{\Omega} \xi(i,\omega) \ g_i(\omega) {\bf P}(d \omega), \eqno(3.1a)
 $$
 $$
  \xi = \vec{\xi} = \{\xi(1), \xi(2), \ldots, \xi(d) \}  =
  \{\xi(1,\omega), \xi(2,\omega), \ldots, \xi(d,\omega) \}
 $$
 is said to be an element of associate space:  $ g = \vec{g} = (g_1,g_2, \ldots, g_d) \in  X'. $\par
 Analogously, arbitrary element of conjugate space $ X^{(d),*}, $ i.e. $ h=h(\xi) $  may be uniquely
represented on the form

$$
h(\xi) = \sum_{i=1}^d h_i(\xi(i)). \eqno(3.1b)
$$

 It is easy to verify that if the space $ X $ has Absolutely Continuous Norm (ACN):

$$
\forall \eta \in X \ \Rightarrow \lim_{ {\bf P}(A) \to 0+ } \int_A |\eta(\omega)| {\bf P}(d \omega) = 0.
$$
then $ X^{(d),'} = X^{(d),*} $ and any linear continuous functional on the space $ (X^{(d)},B) $
may be uniquely represented by the formula (3.1a). \par

{\bf Proposition 3.1.} \par
{\bf A.} The expression (3.1a) represented an element of associate space $ X^{(d),'} $ iff
$ \forall i \ g_i \in  X' $ and

$$
\max_i || g_i||X' \le ||l_g||X^{(d),'} \le \sum_{i=1}^d || g_i||X'. \eqno(3.2a)
$$

{\bf B.} The expression (3.1b) represented an element of conjugate (dual) space $ X^{(d),*} $ iff
$ \forall i \ h_i \in  X^* $ and

$$
\max_i ||h_i||X^* \le ||h||X^{(d),*} \le \sum_{i=1}^d || h_i||X^*. \eqno(3.2b)
$$
 {\bf Proof} is at the same as the proof of proposition  2.1 and may be omitted. \par

 For instance, let $ X = L_p(\Omega) = L(p), $ where $ 1 < p < \infty; $ denote $ q = p/(p-1). $ Then
$ X' = X^* = L_q(\Omega) $ and we deduce by virtue of proposition 3.1 that every linear continuous
functional on the space $ (L(p)^{(d)} $ may be uniquely represented by the formula (3.1a), where
$ g_i \in L(q) $ and

$$
\max_i || g_i||L(q) \le ||l_g||(L(p)^{(d),*})  \le \sum_{i=1}^d || g_i||L(q). \eqno(3.3)
$$

\vspace{3mm}

\section{Multidimensional tail inequalities.}

\vspace{3mm}

 Let $ D $ be arbitrary central-symmetric convex closed bounded set  with non-empty
interior in the space $ R^d. $  We will  denote by  $ D(u), $ where $ u  $ is "great"
numerical parameter: $ u \ge 1 $ its $ u- $ homothetic transformation:

$$
D(u) = \{x, \ x \in R^d, \ x/u \in D \}. \eqno(4.0)
$$

 We intend to obtain in this section the exponential exact as $  u \to \infty $
estimation for the probability

$$
P_{D,\xi}(u) \stackrel{def}{=}  {\bf P}(\xi \notin D(u))  \eqno(4.1)
$$
under assumption that $ \xi \in X^d. $ \par
 For instance, if $ D $ is an unit Euclidean ball in the space $ R^d, $ then

$$
P_{D,\xi}(u) \stackrel{def}{=}  {\bf P}(|\xi|_2 > u).
$$

 We will use the following fact, see, e.g. \cite{Schaefer1}, chapter 4, section 1:
as long as the space $ R^d $ is finite-dimensional, $  D = D^{oo} = M^o, $  where $ M=D^o $
denotes the {\it polar} of the set $ D:$

$$
M=D^o = \cap_{x \in D} \{y, \ y \in R^d, (x,y) \le 1  \}.
$$

\vspace{3mm}

 Let us denote by $ H(\extr(M),\epsilon) =:H_M(\epsilon) $ the entropy of the set $ \extr(M) $
relative the classical Euclidean distance and set $ N_M(\epsilon) = \exp (H_M(\epsilon)). $ \par

\vspace{3mm}

{\bf Theorem 4.1.} Suppose the random vector $ \xi $ belongs to the space $ L_p^{(d)}. $ If the
following integral converges:

$$
I(M,p) \stackrel{def}{=} \int_0^1  N_M^{1/p}(\epsilon) \ d\epsilon < \infty, \eqno(4.2)
$$
then

$$
P_{D,\xi}(u) \le C(I(M,p),d) \ u^{-p}, u \ge 1. \eqno(4.3)
$$

{\bf Proof.}  We conclude by virtue of definition of the set $ M $

$$
P_{D,\xi}(u) =  {\bf P} \left( \sup_{t \in M} (\xi, t)  > u \right). \eqno(4.4)
$$
 Obviously, instead the set $ M $ in the equality (4.4) may be used the set $ \extr(M) $
of all  {\it extremal points} of the set $ M:$

$$
P_{D,\xi}(u) =  {\bf P} \left( \sup_{t \in \extr(M)} (\xi, t)  > u \right). \eqno(4.5)
$$

 Let us consider the (separable, moreover, continuous) random field $ \eta(t) = (\xi,t), $
where $ t \in \extr(M). $ We have by means of definition of  $ || \cdot||(L_p^{(d)},B) $ norm:

$$
\sup_{t \in \extr(M)} |\eta(t)|_p \le ||\xi||L_p^{(d)} < \infty, \eqno(4.6)
$$
and analogously

$$
\forall t,s \in \extr(M) \Rightarrow |\eta(t) - \eta(s)|_p \le ||\xi||L_p^{(d)} \cdot
|t-s|_2. \eqno(4.7)
$$
  The assertion of theorem 4.1 follows immediately from the main result of paper belonging to
G.Pizier \cite{Pizier1}.\par
\vspace{3mm}
{\bf  Remark 4.1.} The case when

$$
\sup_{t \in \extr(M)} ||\eta(t)||G^{(d)}\psi < \infty
$$
and the correspondent distance

$$
\rho(t, s) = ||\eta(t) - \eta(s) ||G^{(d)}\psi
$$
may be investigated analogously. See detail description  with constant estimates in
\cite{Ostrovsky2},  chapter 3, section 3.17.\par

\vspace{3mm}

{\bf  Remark 4.2.} We conclude as long as $ \extr(M) \subset \partial{D}  $ that if the set $ D $
has a smooth boundary, for instance, of the class $ C^1  $ piece-wise,

$$
N(\epsilon) \le C \cdot \epsilon^{-(d-1)}, \ d \ge 2.
$$
\vspace{3mm}

{\bf  Remark 4.3.} Note that the condition (4.2) is satisfied for all the values $ p \in (0,\infty) $
if for example the set $ \extr(M) $ is finite: $ N_M(\epsilon) \le \card(\extr(M)). $ \par
  This occurs, e.g., when   the set $ M $ is (multidimensional,
in general case) polytop.\par
\vspace{3mm}
{\bf Example 4.1.} Let $ D = B $ be the standard unit ball in the space $ R^d; $ then $ M = \extr(M) =
\partial B $ is unit sphere in this space, the entropy
integral $ I(M,p) $ (4.2) converges iff $ p > d-1 $ and in this case

$$
P_{B,\xi}(u) \le C(I(M,p),d) \ u^{-p}, u \ge 1. \eqno(4.8)
$$
 Evidently, the estimate (4.8) is true for arbitrary convex bounded domain $ D $ such that

$$
\sup_{x \in D} |x| \le 1.
$$

\vspace{3mm}

{\bf  Remark 4.3.} Assume that the condition (4.2) is satisfied for any {\it diapason } $  (a,b) $
of a values $ p; $ here $ d-1 < a < b \le \infty; $  then

$$
P_{D,\xi}(u) \le \inf_{ p \in (a,b) } \left[ C(I(M,p),d) \ u^{-p} \right], u \ge 1. \eqno(4.9)
$$
 When $ b = \infty, $  then it may be obtained from (4.9) the exponential decreasing as $ u \to \infty $
estimate for the probability  $ P_{D,\xi}(u). $ \par

\vspace{3mm}

\section{Moment inequalities for sums of random vectors.}

\vspace{3mm}

 We recall before formulating the main result  some useful for us
moment inequalities  for the  sums of centered  martingale differences (m.d.) $ \theta(i) $
relative some filtration $ \{ F(i) \}: $

$$
F(0) = \{ \emptyset, \Omega \}, F(i) \subset F(i+1) \subset F:
$$

 $  \forall k=0,1,\ldots, i-1 \ \Rightarrow $
$$
{\bf E}\theta(i)/F(k)=0; \ {\bf E} \theta(i)/F(i)= \xi(i) \ (\mod {\bf P}),
$$
and for the independent r.v., \cite{Ostrovsky4}. Namely, let $ \{ \theta(i) \} $ be a
sequence of centered  martingale differences relative any filtration; then

$$
\sup_n \sup_{b \in S(1)} \left|\sum_{i=1}^n b(i) \theta(i) \right|_p \le K_M(p) \
\sup_i |\theta(i)|_p, \eqno(5.0)
$$
where for the {\it optimal value} of the constant $ K_M = K_M(p) $ there holds  the inequality

$$
 K_M(p) \le   p \ \sqrt{2}, \ p \ge 2.
$$
 Note that the upper bound in (5.0)
 $$
  K_I(p) \le  0.87 p/\log p, \ p \ge 2
 $$
is true for the independent centered r.v. $ \{\theta(i)) \}, $ see also \cite{Ostrovsky4}. \par
 Applying the inequality (5.0) for the value  $ n=d,$ we obtain the following result. \par
\vspace{3mm}
{\bf  Proposition 5.1.} \par
{\bf A.} Suppose the coordinates of the vector
$ \xi = \vec{\xi} = \{\xi(1), \xi(2), \ldots, \xi(d) \} $ are centered martingale differences. Then

$$
||\vec{\xi}||(L(p)^{(d)}) \le  K_M(p) \max_i |\xi(i)|_p. \eqno(5.1a)
$$

{\bf B.} Suppose the coordinates of the vector
$ \vec{\xi} = \vec{\xi} = \{\xi(1), \xi(2), \ldots, \xi(d) \} $ are centered independent r.v. Then

$$
||\vec{\xi}||(L(p)^{(d)}) \le  K_I(p) \max_i |\xi(i)|_p. \eqno(5.1b)
$$
\vspace{3mm}
 Let us denote for any function $ \psi(\cdot) \in G\Psi $

 $$
 \psi_K(p) = K_M(p) \psi(p), \ \psi_I(p) = K_I(p) \psi(p).
 $$
A small consequence of a proposition 5.1: \par
{\bf  Proposition 5.2.} \par

{\bf A.} Suppose the coordinates of the vector
$ \vec{\xi} = \xi = \{\xi(1), \xi(2), \ldots, \xi(d) \} $ are centered martingale differences. Then

$$
||\vec{\xi}||(G\psi_M^{(d)}) \le  \max_i ||\xi(i)||G\psi. \eqno(5.2a)
$$

{\bf B.} Suppose the coordinates of the vector
$ \vec{\xi} = \vec{\xi} = \{\xi(1), \xi(2), \ldots, \xi(d) \} $ are centered independent r.v. Then

$$
||\vec{\xi}||(G\psi_I^{(d)}) \le  \max_i ||\xi(i)||G\psi. \eqno(5.2b)
$$

 We intend now to generalize the famous Rozenthal's inequality on the multidimensional case.
 Let again $ \psi \in \Psi  $ and let
 $ \eta = \eta(1) $ be a centered random vector from the space $ G\psi^{(d)};  $ let
 $  \eta(2), \eta(3), \ldots, \eta(n) $ be independent copies $ \eta. $ We denote

 $$
 \zeta(n) = n^{-1/2} \sum_{j=1}^n \eta(j).
 $$
{\bf  Proposition 5.3.} \par

$$
\sup_n ||\zeta(n)||(G\psi_I^{(d)}) \le ||\eta||(G\psi^{(d)}). \eqno(5.3)
$$
{\bf Proof} follows immediately from the classical Rozenthal's inequality, see, e.g.
\cite{Johnson1}, \cite{Johnson2},  \cite{Ostrovsky4},  \cite{Rozenthal1},
\cite{Sharachmedov1}, \cite{Utev1}. Namely, let $ b $ be arbitrary deterministic vector
from the set $ S(d). $  We use the Rozenthal's inequality for the one-dimensional mean zero
independent r.v. $ \nu(j) = (\eta(j),b): $

$$
\sup_n | n^{-1/2} \sum_{j=1}^n \nu(j)|_p \le K_I(p) \cdot |\nu|_p. \eqno(5.4)
$$
 The assertion (5.3) follows from (5.4) after dividing on the $ \psi(p) $ and
taking maximum over $ b; \ b \in S(d) $ and $ p; \psi(p) \in (0,\infty). $ \par
Note that the estimate of a view

$$
\sup_n ||\zeta(n)||(G\psi_M^{(d)}) \le ||\eta||(G\psi^{(d)}) \eqno(5.5)
$$
is true for the sequence of (centered) martingale differences $ \{ \eta(j)\}, j=1,2,\ldots,n. $ \par

\vspace{3mm}

\section{ Fundamental function for m.d.r.i. spaces. }

\vspace{3mm}

 Recall that the  fundamental function $ \chi_X = \chi_X(\delta), \ \delta \in (0, \mu(\Omega)) $
for the r.i. space with measure $ \mu(\cdot) \ (X, ||\cdot||X) $  is defined by the formula

$$
 \chi_X(\delta) = \sup_{A, {\bf P}(A) \le \delta} ||I(A)||X, \eqno(6.0)
$$
where as ordinary $ I(A) $ is the indicator function of the event $ A. $ \par
 We intend to generalize this definition in the multidimensional case of the space $ X^{(d)}. $
Namely, let $ \vec{A} = \{ A(i) \}, i=1,2,\ldots,d $  be a {\it family } of measurable subsets of
the whole space $ \Omega $  and define the vector-function

$$
I(\vec{A}) = \{I(A(1), I(A(2)), \ldots, I(A(d)) \}.
$$

 The fundamental function of the m.d.r.i. space $ X^{(d)} \
\chi_{X^{(d)}}(\delta_1,\delta_2, \ldots, \delta_d), \ \delta_i \in [0,1] $
may be defined as follows:

$$
\chi_{X^{(d)}}(\delta_1,\delta_2, \ldots, \delta_d) =
\sup_{A(i): {\bf P}(A(i)) \le \delta_i} ||I(\vec{A})||(X^{(d)})=
$$
$$
\sup_{A(i): {\bf P}(A(i)) \le \delta_i} \sup_{b \in S(d)} ||\sum_{i=1}^d b(i) I(A(i))||X. \eqno(6.1.)
$$

{\bf Proposition 6.1.}

$$
\max_i \chi_X(\delta_i) \le  \chi_{X^{(d)}}(\delta_1,\delta_2, \ldots, \delta_d) \le
\sum_{i=1}^d \chi_X(\delta_i). \eqno(6.2)
$$

 {\bf Remark 6.1.} Note that the lower bound in the bilateral inequality (6.2) is attained, for
instance, when $ X = L_p(\Omega), \ p \ge 2 $   and when the sets $  A(i) $ are disjoined. \par

\vspace{3mm}

\section{ Multidimensional spaces of subgaussian  and pre-gaussian type. }

\vspace{3mm}

 We consider in this section the case when  the space $ (X, ||\cdot||X)  $ coincides with
the space $ \Phi(\phi) $ for some $ \phi \in \Phi. $ \par
 Let the random vector $ \vec{\xi} = \xi $ belongs to the space $ \Phi^{(d)}. $ Recall that
this imply that $ {\bf E} \vec{\xi} = {\bf E}\xi = 0 $ and $ \forall \lambda \in R $

$$
\sup_{b \in S(d)} {\bf E} \exp \{\lambda ( \xi,b) \} \le
\exp \{ \phi(\lambda \cdot ||\xi||\Phi^{(d)}) \},
$$
or equally

$$
\forall \mu \in R^d \ \Rightarrow
{\bf E} \exp \{(\xi,\mu) \} \le  \exp\{\phi(|\mu|_2 \cdot ||\xi||\Phi^{(d)}) \}. \eqno(7.1)
$$
 Such a random vectors are called pre-gaussian. \par
 We intend now to generalize the last definition. Let $ D $ be positive definite symmetrical constant
matrix (linear operator) of a size $ d \times d. $   By definition, the (necessary mean zero)
random vector $ \xi = \vec{\xi} $ belongs to the space $ \Phi_D(\phi), \ \phi \in \Phi, $ if there
is a non-negative constant $  \tau  $ which dependent only on the distribution of the r.v.
$ \xi:  \tau = \tau(\Law(\xi)) $ such that $ \forall \lambda \in R $ and for arbitrary  $ b \in S(d) $

$$
{\bf E} \exp (\lambda (\xi,b)) \le \exp \phi \left( \lambda \tau \sqrt{ (Db,b)} \right). \eqno(7.2)
$$

{\bf Definition 7.1.} \par
 The minimal value of the constant $ \tau, \ \tau \ge 0 $ for which the inequality 7.2 is satisfied for all
the prescribed values $ \lambda, b $  is said to be  the $ \phi,D $ norm of the random vector $ \xi: $

$$
||\xi||_{\phi,D} = \inf \{\tau, \ \tau > 0,  \
 \forall \lambda \in R, \ \forall b \in S(d) \ \Rightarrow
$$

$$
{\bf E} \exp (\lambda (\xi,b)) \le \exp \phi \left( \lambda \tau \sqrt{ (Db,b)} \right) \} \eqno(7.3)
$$
or equally

$$
||\xi||_{\phi,D} = \inf_{\lambda \ne 0, \ b \in S(d)}
\frac{\phi^{-1} (\log {\bf E} \exp(\lambda (\xi,b)))}{|\lambda| ||b||_D},  \eqno(7.4)
$$
here and hereafter we denote

$$
||b||_D = \sqrt{(Db,b)}.
$$
 We will denote the Banach space of all such a random vectors (i.e. with finite norm $ ||\xi||_{\phi,D} )$
as $ \Phi^{(d)}(\phi,D). $ \par
 For instance, if a random vector $ \eta $ has Gaussian centered distribution with variation
 $ R: \Law(\xi) = N(0,R), $ then we can take $ D=R $ and $ \phi(\lambda):=
 \phi_0(\lambda) \stackrel{def}{=} 0.5 \lambda^2: $

$$
||\eta||_{\phi_0, R} = 1.
$$
{\bf Definition  7.2.} \par
 The random  vector $  \xi $ is said to be subgaussian relative the symmetric positive definite
matrix $  D  $, if in the inequality (7.3) it can be taken  $ \phi(\lambda) = 0.5 \lambda^2: $

$$
{\bf E} \exp (\lambda (\xi,b)) \le \exp  \left( 0.5 \lambda^2 \tau^2 (Db,b) \right) \} \eqno(7.5)
$$
and the random  vector $  \xi $ is said to be strong subgaussian,
if in the inequality (7.5) it can be taken  $ \phi(\lambda) = 0.5 \lambda^2 $ and $ D = R = \Var(\xi): $

$$
{\bf E} \exp (\xi,\mu) \le \exp (0.5 (R \mu,\mu) ). \eqno(7.6)
$$
 This definitions belong to V.V.Buldygin and Yu.V.Kozatchenko, see \cite{Buldygin1}, \cite{Buldygin2},
where are described some applications. Another investigations and applications  see in
\cite{Ostrovsky7}. \par

{\bf Definition  7.3.} \par
Let $ \phi $ be any function from the set $  \Phi. $
 The random  vector $  \xi $ is said to be $ \phi- $ strong subgaussian,
if in the inequality (7.3) it can be taken  $ D = R = \Var(\xi): $

$$
{\bf E} \exp (\lambda (\xi,b)) \le \exp \phi \left( \lambda \tau \sqrt{ (Rb,b)} \right).  \eqno(7.7)
$$
\vspace{4mm}

  Too more general definition. Let $ \nu = \nu(\mu), \mu \in R^d $ be
even: $ \nu(-\mu) = \nu(\mu) $ strong convex which takes positive values for non-zero
 arguments twice continuous differentiable function, such that
$$
 \nu(0) = 0, \ \grad \nu(0)  = 0, \ \lim_{r \to \infty}\min_{|\mu|_2 \ge r} ||\grad \nu(\mu)||_2 = \infty;
 \eqno(A)
$$

$$
\inf_{\mu \in R^d} \lambda_{\min} \left[ \frac{\partial^2 \nu}{\partial \mu_j \partial \mu_k} \right] > 0;
\eqno(B)
$$

$$
 \nu(\mu(k)\times e(k)) = \nu(0,0, \ldots,0,\mu(k),0,0, \ldots, 0)
  \le \nu(\vec{\mu}), \eqno(C)
$$
where

$$
e(k) \stackrel{def}{=} (0,0, \ldots,0,1,0,0, \ldots, 0),
$$
and $ "1" $ stands on the place $ k. $ \par
 We denote the set of all such a function as $ \Phi^{(d)}; \ \Phi^{(d)} =\{ \nu(\cdot) \}. $ \par
 For instance, the  conditions (A), (B) and (C) are satisfied for non-degenerate centered (multidimensional)
Gaussian distribution. \par
{\bf  Definition 7.4.} \par
 We will say that the {\it centered} random vector $ \eta = \eta(\omega) $
belongs to the space $ \Phi^{(d)}(\nu), $ if there exists some non-negative constant
$ \tau \ge 0 $ such that

$$
\forall \mu \in R^d  \ \Rightarrow
{\bf E} \exp(\mu, \eta) \le \exp[ \nu(\mu \ \tau) ].
$$
 The minimal value $ \tau $ satisfying this inequality  is called a $ \Phi^{(d)}(\nu) \ $ norm
of the variable $ \eta, $ write
 $$
 ||\eta||\Phi^{(d)}(\nu) = \inf \{ \tau, \ \tau > 0: \ \forall \mu \in R^d \ \Rightarrow
 {\bf E}\exp(\mu \eta) \le \exp(\nu(\mu \ \tau)) \}.
 $$

{\bf Theorem 7.1.} The set $ \Phi^{(d)}(\nu) $ with ordinary operation  equipped with correspondent
norm  $ ||\cdot||\Phi^{(d)}(\nu) $ is (complete) r.i. Banach  space over $ (\Omega,F,{\bf P} ). $\par
{\bf Proof } is at the same as in one-dimensional case, see \cite{Kozatchenko1}, \cite{Ostrovsky2},
 \ chapter 1; see also  \cite{Bagdasarova1}, and may be omitted. \par

\vspace{3mm}
  Recall that the {\it multidimensional } Young-Fenchel, or Legendre transform $ \nu^*(x), \ x \in R^d $
of a function $ \nu: R^d \to R  $  is defined by the formula

$$
\nu^*(x) = \sup_{\mu \in R^d} (x \mu - \nu(\mu)). \eqno(7.9)
$$
  It is well known, \cite{Ioffe1}, chapter 5, that $ \nu^*(x) $ is convex function and if the function
$ \phi(\cdot) $  is convex and  continuous on the its support, then

$$
\nu^{**}(\mu) = \nu(\mu), \ \mu: \nu(\mu) < \infty
$$
(theorem of Fenchel-Moraux.) \par
 The multidimensional tail function for the random vector $  \xi \ T_{\xi}(x), x \ge 0, x \in R^d $
may be defined as follows:

$$
T_{\xi}(x) = \max_{\pm}{\bf P}(\pm \xi(1)| \ge x_1,  \pm \xi(2) \ge x_2, \ldots, \pm \xi(d) \ge x_d),
$$
where the exterior maximum is calculated over all combinations the signs $ \pm. $ \par
 {\bf Theorem 7.2.} The non-zero centered random vector $ \xi $ belongs to the space
$ \Phi^{(d)}(\nu,D) $ iff

$$
T_{\xi}(\vec{x}) \le \exp \left( - \nu^*(C \vec{x}) \right) \eqno(7.10)
$$
(Tchernoff-Tchebychev's  inequality). \par
 {\bf Proof.} \par
{\sc A.}  Let $ \xi \in \Phi^{(d)}(\nu) $ and  $ ||\xi||\Phi^{(d)}(\nu) = 1. $
Denote as before also $ x = \vec{x} = \{x(1), x(2), \ldots, x(d) \}, $ and
assume without loss of generality that all the coordinates of the vector $ x $ are strictly
positive. \par
 As long as in the case $ \max(x(i)) < 1 $ the inequality (7.10) is evident, we consider further
a possibility $ \max(x(i)) \ge 1. $ \par

  We have by virtue of Tchebychev's inequality:

$$
T_{\xi}(x) \le {\bf E} \exp(\mu,\xi) / \exp(\mu,x) \le \exp \left(-(\mu,x) - \nu(\mu)) \right). \eqno(7.11)
$$
 The assertion (7.10) is obtained from (7.11) after minimization over $ \mu; $
in the considered case it may be adopted  $ C = 1.$ \par

{\sc B.} Conversely, let the estimate (7.10) there holds with unit constant $ C: C = 1. $
Assume for definiteness $ \vec{\mu}\ge 0. $
As before, it is sufficient to consider only the case $ \forall j=1,2,\ldots,d \ \mu(j) \ge 1. $\par
 We conclude after integration by parts:

$$
{\bf E} \exp(\mu, \xi) \le C + \prod_{i=1}^d \mu(i) \cdot \int_{R^d_+}
\exp( (\mu,x) - \nu^*(x) ) \ dx. \eqno(7.12)
$$
 We can estimate the last integral by means of saddle-point method,
 \cite{Fedoruk1}, chapter 2; see also  \cite{Kozatchenko1}:

$$
{\bf E} \exp(\mu \xi) \le C + \prod_{i=1}^d \mu(i) \cdot
\sup_x \exp(C_2 (\mu,x) - \nu^*(x) ) =
$$

$$
 C + \prod_{i=1}^d \mu(i) \cdot \nu^{**}(C_3 \mu) \le \nu^{**}(C_4 \mu) =
\nu(C_4 \mu);
$$
we used the theorem of Fenchel-Moraux. \par
 This completes the proof of theorem 7.2. \par

 {\bf Example 7.1.} \par

{\bf A. } Suppose for some $ \tau = \const > 0, \ \phi \in \Phi $ and for any  strictly
positive definite symmetrical matrix $  D $  of a size $ d \times d $

$$
\forall \mu \in R^d \Rightarrow {\bf E } \exp(\xi, \mu) \le \exp(\phi(\tau \cdot |\mu|_D )),
$$
then

$$
T_{\xi}(x) \le  \exp(-\phi^*(x \cdot |\mu|_{D^{-1}}/\tau )).
$$

{\bf B.} Conversely, if $ {\bf E} \xi = 0 $ and  for some constant $ \tau > 0 $

$$
T_{\xi}(x) \le  \exp(-\phi^*( x \cdot |\mu|_{D^{-1}}/\tau )),
$$
then there is a positive finite constant $ C = C(\phi) $ such that

$$
\forall \mu \in R^d \Rightarrow {\bf E } \exp(\xi, \mu) \le \exp(\phi(C \tau \cdot |\mu|_D )).
$$
 In order to show the precision of result of theorem 7.2, let us consider the following \\
{\bf Example 7.2.}  Let  $ (\xi(1), \xi(2)) $ be two-dimensional mean zero Gaussian distributed
normed random vector such that

$$
{\bf E} \xi^2(1) = {\bf E} \xi^2(2)=1, \ {\bf E} \xi(1)\xi(2) =:\rho \in (-1,1).
$$
 It follows from the theorem 7.2 that at $ x_1 > 0, \ x_2 > 0 $

$$
{\bf P} (\xi(1) > x_1,\xi(2) > x_2) \le
\exp \left(-0.5 (1-\rho^2)^{-1}(x^2_1 - 2 \rho x_1 x_2 + x_2^2 ) \right),
$$
but as $ x_1 \to \infty, x_2 \to \infty $ independently inside the fixed angle

$$
\arctan \rho + \varepsilon_0 \le \frac{x_2}{x_1} \le \arctan 1/\rho - \varepsilon_0,
$$

$  \varepsilon_0 = \const \le 0.25( \arctan 1/\rho - \arctan \rho ):  $

$$
{\bf P} (\xi(1) > x_1, \xi(2) > x_2) \sim C_{1,2}(\rho) \ x_1^{-1} x_2^{-1} \
\exp \left(-0.5 (1-\rho^2)^{-1}(x^2_1 - 2 \rho x_1 x_2 + x_2^2 ) \right).
$$

\vspace{4mm}
 We continue the investigation of these spaces.
Let the random vector $ \vec{\xi} = \xi $ belongs to the space $ \Phi^{(d)}(\nu) $ and
let $ k=1,2,\ldots,d. $ We denote

$$
\nu_k(z) = \nu(z \times e(k)), \ z \in R,
$$
so that

$$
{\bf E} \exp(z \xi(k)) \le \exp( \nu_k(z) );
$$

$$
\psi_k(p) := \frac{p}{\nu_k^{-1}(p)}.
$$
\vspace{3mm}
 {\bf Theorem 7.3.} Suppose the function $ \nu(\cdot)  $  satisfies the conditions (A), (B) and (C).
On the (sub-)space of all {\it centered} random vectors  $ \{ \xi = \vec{ \xi} \} $
the norm  $ ||\xi|| \Phi^{(d)}(\nu) $  is equivalent to the other, so-called "moment norm":

 $$
|||\xi|||G\vec{\nu} \stackrel{def}{=} \sum_{k=1}^d   \sup_{p \ge 1} \frac{|\xi(k)|_p}{\psi_k(p)}:
\eqno(7.13)
 $$

 $$
C_1(\nu) ||\xi|| \Phi^{(d)}(\nu) \le ||\xi||G\vec{\nu} \le C_2(\nu)||\xi|| \Phi^{(d)}(\nu).
\eqno(7.14)
 $$
{\bf Proof.} \par

{\sc A.} Let  $ ||\xi||\Phi^{(d)}(\nu) = 1. $  From the direct definition of the norm
$ ||\xi|| \Phi^{(d)}(\nu) $ it follows

$$
\forall \mu \in R^d  \ \Rightarrow {\bf E} \exp(\mu, \xi) \le \exp[ \nu(\mu) ].\eqno(7.15)
$$
  We substitute into inequality the value $ \mu = z\cdot e_k, \ z \in R:$

$$
{\bf E} \exp(z \xi(k)) \le \exp( \nu_k(z) ). \eqno(7.16)
$$
 The inequality (7.16) means that the one-dimensional  r.v. $ \xi(k) $ belongs to the space
$ \Phi(\nu_k) $  and has there unit norm. \par
 From the theory of these spaces  \cite{Kozatchenko1}, \cite{Ostrovsky2}, chapter 1 it follows

 $$
\sup_{p \ge 1} \frac{|\xi(k)|_p}{\psi_k(p)} =: Y(k) < \infty,
 $$
therefore
$$
|||\xi|||G\vec{\nu} \le \sum_{k=1}^d Y(k) < \infty,
$$
since the last sum is finite.\par

{\sc B.} Let now $ {\bf E }\xi = 0 $ and $ |||\xi|||G\vec{\nu} < \infty; $ without loss of
generality we can and will suppose $ |||\xi|||G\vec{\nu} =1. $ \par
 It follows from the last inequality that

 $$
\sup_{p \ge 1} \frac{|\xi(k)|_p}{\psi_k(p)} \le 1
 $$
and hence

$$
||\xi(k)||\Phi(\nu_k) =:Z_k < \infty.
$$
 We conclude  by virtue of condition (C) of our theorem that

$$
|| \xi(k) \cdot e(k)|| \Phi^{(d)}(\nu) = ||\xi(k)||\Phi(\nu_k) = Z_k.
$$
 As long as $ \xi = \sum_{k=1}^d \xi(k) \cdot e(k), $ we obtain using triangle inequality

$$
|| \xi|| \Phi^{(d)}(\nu) \le \sum_{k=1}^d ||\xi(k)\times e(k)||\Phi^{(d)}(\nu) =
$$

$$
 \sum_{k=1}^d ||\xi(k)||\Phi(\nu_k) = \sum_{k=1}^d Z_k < \infty,
$$
Q.E.D. \par

\vspace{4mm}
 {\bf Proposition 7.1.} Suppose  $ \xi \in \Phi^{(d)}(\nu), \xi \ne 0. $   Denote

 $$
 \overline{\nu}(\mu) = \sup_{n=1,2,\ldots} n \ \nu(\mu/\sqrt{n}).
 $$
 Note that the function $ \overline{\nu}(\mu) $ there exists and satisfies the conditions
 (A), (B) and (C) since there exists a limit

 $$
 \lim_{n \to \infty} n \ \nu(\mu/\sqrt{n})= 0.5 \ (R\mu, \mu), \ R = \Var{\xi}.
 $$
 Let $ \{ \xi^{(j)} \},  \ j=1,2,\ldots $ be independent copies of $ \xi. $
Let us denote

$$
S(n) =  n^{-1/2} \sum_{j=1}^n  \xi^{(j)}.
$$

We observe as in the one-dimensional case:

$$
\sup_n || S(n) ||\Phi^{(d)}({\overline{\nu}}) \le ||\xi||\Phi^{(d)}(\mu).
\eqno(7.17)
$$
 As a consequence we obtain an exponential bounds for normed sums of independent  centered
random vectors.  \par

{\bf Proposition 7.2.}
  Let $ \{ \xi^{(j)} \},  \ j=1,2,\ldots $ be independent copies of mean zero non-trivial
random vector  $ \xi $ belonging to the space $ \Phi^{d)}. $
We have the following uniform tail estimate for the norming sum of the r.v.
$ \{ \xi^{(j)} \}:$

$$
\sup_n T_{S(n)}( \vec{x}) \le \exp
\left( - \overline{\nu}^* \left( \vec{x}/||\xi||\Phi^{(d)}(\nu) \right) \right), \ \vec{x} > 0. \eqno(7.18)
$$

\vspace{3mm}

\section{ Multidimensional random fields (processes). }

\vspace{3mm}
 Let $ Y = \{ y \} $ be arbitrary set and let (in this section) $ \xi = \vec{\xi} = \vec{\xi}(y)  $
be separable random field (process) with values in the space $ R^d: $

$$
\vec{\xi}(y) = \{ \xi(1,y), \xi(2,y), \ldots, \xi(d,y) \}.
$$

 The aim of this section is estimate the joint distribution (more exactly, joint tail function)
as $ \vec{v} = \{ v_1,v_2, \ldots, v_d \} \to \infty \Leftrightarrow \min v_i \to \infty $
for the coordinate-wise maximums:

$$
U(\vec{v}) \stackrel{def}{=} {\bf P}
\left( \sup_{y \in Y}\xi(1,y) > v_1, \sup_{y \in Y}\xi(2,y) > v_2, \ldots,
 \sup_{y \in Y}\xi(d,y) > v_d \right) \eqno(8.0)
$$
or for brevity

$$
U(\vec{v}) = {\bf P}
\left( \sup_{y \in Y} \vec{\xi}(y) > \vec{v}\right),  \vec{v} = (v_1,v_2, \ldots,v_d),
$$
$$
\sup_{y \in Y} \vec{\xi}(y) = \{\sup_{y \in Y}\xi(1,y), \sup_{y \in Y}\xi(2,y),
\ldots, \sup_{y \in Y}\xi(d,y) \}.
$$
 We denote

$$
\overline{\xi}(j) = \sup_{y \in Y}\xi(j,y),
$$
so that

$$
U(\vec{v}) \stackrel{def}{=} {\bf P}
\left( \overline{\xi}(1) > v_1,  \overline{\xi}(2) > v_2, \ldots,   \overline{\xi}(d) > v_d \right).
$$

 The one-dimensional case $ d = 1 $  with described applications is in detail investigated in
\cite{Ostrovsky2}, introduction and chapter 4. \par
 Notice that in a particular case when  $ v_1 = v_2 = \ldots = v_d = v $ the probability
 $ U(\vec{v}) $ has a form:
 $$
 U(v,v,\ldots,v) = {\bf P}( \min_{i=1,2,\ldots,d} \sup_{y \in Y} \xi(i,y) > v),
 $$
i.e. the probability  $ U(\vec{v}) $  is the tail distribution for {\it minimax } of the random
field $ \xi(i,y). $ \par
 We assume that

 $$
 \forall y \in Y \ \Rightarrow \vec{\xi}(y) \in \Phi^{(d)}(\nu)
 $$
and moreover

$$
\sup_{ y \in Y} ||\vec{\xi}(y)||\Phi^{(d)}(\nu) = 1. \eqno(8.1)
$$

 Let us introduce the following {\it natural } distance, more exactly, {\it bounded:}
$ d(y_1,y_2) \le 2 $ semi-distance on the set $  Y: $

$$
d(y_1,y_2) = ||\vec{\xi}(y_1) - \vec{\xi}(y_2)||\Phi^{(d)}(\nu). \eqno(8.2)
$$
 Denote as usually for any subset $ V, \ V \subset Y $ the so-called
{\it entropy } $ H(V, d, \epsilon) = H(V, \epsilon) $ as a logarithm
of a minimal quantity $ N(V,d, \epsilon) = N(V,\epsilon) = N $
of a (closed) balls $ S(V, t, \epsilon), \ t \in V: $
$$
S(V, t, \epsilon) \stackrel{def}{=} \{s, s \in V, \ d(s,t) \le \epsilon \},
$$
which cover the set $ V: $
$$
N = \min \{M: \exists \{t_i \}, i = 1,2,…, M, \ t_i \in V, \ V
\subset \cup_{i=1}^M S(V, t_i, \epsilon ) \},
$$
and we denote also for the values $ \epsilon \in (0,1) $
$$
H(V,d,\epsilon) = \log N; \ S(t_0,\epsilon) \stackrel{def}{=}
 S(Y, t_0, \epsilon), \ H(\epsilon) = H(d, \epsilon) \stackrel{def}{=} H(Y,d,\epsilon).
$$
 Here $ t_0 $ is the so-called "center" of the set $ Y $ relative the distance $ d, $
i.e. the point for which

$$
\sup_{t \in Y} d(t_0, t) \le 1.
$$
 We can and will assume the completeness of the space $ Y $  relative the distance $ d; $
the existence of the center $ t_0 $ it follows from the Egoroff's theorem.\par
 It follows from Hausdorff's theorem that
$ \forall \epsilon > 0 \ \Rightarrow H(Y,d,\epsilon)< \infty $ if and only if the
metric space $ (Y, d) $ is pre-compact set, i.e. is the bounded set with compact closure.\par
 Let $ p $ be arbitrary number from the interval $ (0,1/2). $ Define a function

$$
w(p) = (1-p)\sum_{n=1}^{\infty}  p^{n-1}H(p^n). \eqno(8.3)
$$
 We assume that
 $$
 \exists p_0 \in (0,1/2) \ \forall p \in (0,p_0) \ \Rightarrow w(p) < \infty,
 \eqno(8.4)
 $$
(The so-called "entropy condition"). \par

 The condition (8.4) is satisfied if for example
 $$
 H(\epsilon) \le C + \kappa |\log \epsilon|;  \eqno(8.5)
 $$
 in this case

$$
w(p) \le C + \frac{\kappa |\log p|}{1-p}.
$$
 The minimal value $ \kappa $ for which the inequality (8.5) holds (if there exists)
is said to be {\it entropy dimension } the set $ Y $ relative the distance $ d, $ write:

$$
\kappa = \edim(Y;d).
$$
 For instance, if $ Y $ is closed bounded subset of whole space $ R^m $  with non-empty
interior and for some positive finite constant $ C_1, C_2; \ 0< C_1 \le C_2 < \infty $

$$
 C_1 |y_1 - y_2|^{\alpha}  \le d(y_1,y_2) \le   C_2 |y_1 - y_2|^{\alpha}, \ \alpha = \const
 \in (0,1],
$$
where $ |z| $ is ordinary Euclidean norm, then $  \kappa = \edim(Y;d) = d/\alpha. $ \par

{\bf Theorem 8.1.} Suppose the function $ \nu(\cdot) $ satisfies the conditions (A), (B), (C).
Let also the entropy condition (8.4) be satisfied.  Then the vector random  field $ \vec{\xi}(y) $
is $ d -$ continuous with probability one:

$$
{\bf P}( \vec{\xi}(\cdot) \in C(Y,d) = 1)
$$
and moreover
$$
U(v) \le \inf_{p \in (0, p_0) } \exp \left[w(p) - \nu^*( (v(1-p))) \right]. \eqno(8.6)
$$

{\bf Proof.}  The  $  d - $ continuity of  r.f. $ \xi(\cdot) $ follows immediately from
the main result of the paper \cite{Ostrovsky7}; it remains to prove the estimate (8.6).\par
 We denote as $ S(\epsilon) $ the minimal $ \epsilon - $ net of the set $  Y $
relative the distance $ d: \ \card(S(\epsilon)) = N(\epsilon); $ not necessary to be unique,
but non-random. By definition, $ S(1) = \{t_0\} $ and $ S_n = S(p^n) $ for arbitrary values
$ p $ inside the interval $ (0,p_0). $   \par
 We define for any $ y \in Y $ and $ n=0,1,2,\ldots $ the  non-random functions $ \theta_n(\cdot), $
(projection into the set $ S_n) $ also not necessary to be unique, as follows: $ \theta_0(y) = y_0 $
and
 $$
 \theta_n(y)  = y_j, \ y_j \in S_n, \ d(y,y_j) \le p^n, \ n=1,2,\ldots. \eqno(8.7)
 $$
 The set

$$
S  = \cup_{n=0}^{\infty} S(p^n) \eqno(8.8)
$$
is enumerate dense subset of the space $ (V,d); $ we can choose the set $  S $ as the set of
separability for the random field $ \vec{\xi}(t). $  We have:

$$
\sup_{y \in Y} \xi(y) = \lim_{n \to \infty} \max_{y \in S_n} \xi(y).
$$
 Further,  when $ n \ge 1 $

$$
 \max_{y \in S_n} \xi(y) \le \max_{y \in S_n} (\xi(y) - \xi(\theta_{n-1} y) ) +
 \max_{y \in S_{n-1}} \xi(y),
$$
therefore

$$
\sup_{y \in Y} \xi(y) \le \sum_{n=0}^{\infty} \eta_n, \eqno(8.9)
$$
where

$$
\eta_0 = \xi(y_0), \ \eta_n = \max_{y \in S_n} ( \xi(y) - \xi(\theta_{n-1}y)), \ n=1,2,\ldots.
$$
 We have for positive values $ \mu = \vec{\mu}: \
 {\bf E} \exp(\mu \eta_0) \le \exp(\nu(\mu) ),\ {\bf E} \exp(\mu \eta_n) \le  $

$$
 \sum_{y \in S_n } {\bf E} \exp \left( (\mu,\xi(y) - \xi(\theta_{n-1} y))  \right) \le
 \sum_{y \in S_n} \exp(\nu(\mu \cdot p^{n-1}) ) \le
$$
$$
 N(p^n) \exp \left(\nu(\mu \cdot p^{n-1}) \right) =
 \exp \left( H(p^n) + \nu(\mu \cdot p^{n-1})   \right).
$$

 As long as
$$
{\bf E} \exp( \mu \cdot \sup_{y \in Y} \xi(y) ) \le {\bf E} \exp(\mu \cdot \sum_n \eta_n)=
{\bf E}\prod_{n=0}^{\infty}  \exp(\mu \cdot \eta_n).
$$
We obtain using H\"older's inequality in which we choose $ 1/r(n) =p^{n}(1-p), $ so that
$ r(n) > 1, \ \sum_{n=0}^{\infty}1/r(n) = 1: $

$$
{\bf E} \exp( \mu \cdot \sup_{y \in Y} \xi(y) ) \le
 \prod_{n=0}^{\infty} \left[ {\bf E} \exp( \mu \cdot r(n) \cdot \eta_n  ) \right]^{1/r(n)} \le
$$

$$
\prod_{n=0}^{\infty} \left[\exp( H(p^n) + \nu(r(n) \cdot \mu \cdot p^{n}) )   \right]^{1/r(n)}
\le
$$

$$
\exp \left\{w(p) +   \nu \left( \frac{\mu}{1-p}  \right)  \right\}. \eqno(8.10)
$$

  It remains to use Tchebychev's - Tchernoff's  inequalities in order to obtain the assertion
(8.6) of theorem 8.1 for arbitrary fixed value $ p $ and with consequent optimization over $ p. $ \par

 {\bf Corollary 8.1.} Let us denote

 $$
 \pi(v) = \frac{C}{(\nabla \nu^*(v), v)}
 $$
and assume that $ \pi(v) $ there exists and

$$
\lim_{\min v(i) \to \infty} \pi(v) = 0.
$$
 We deduce choosing in (8.6) the value $ p = \pi(v) $ for all sufficiently great values
$ v = \vec{v}: \ \pi(v) \le p_0: $

 $$
U(v) \le \exp \left[ w \left(\frac{C_2}{\pi(v)} \right) - \nu^*(v) \right]. \eqno(8.11)
 $$

{\bf  Example 8.1.}   Suppose in addition to the conditions of theorem 8.1 the condition 8.5
holds. Then
\vspace{3mm}
$$
U(v) \le C(\kappa) \ \pi(v)^{-\kappa} \ \exp \left( - \nu^*(v)  \right), \ \min_i v_i > v_0.
\eqno(8.12)
$$

\vspace{3mm}
{\bf  Example 8.2.} Let $ \xi(y) = (\xi(1,y), \xi(2,y)), \ y \in Y $ be a two-dimensional
 centered Gaussian random field with the following covariation symmetrical positive definite
 matrix-function: $  R = R(y_1,y_2) = \{ R_{i,j}(y_1,y_2) \},  i,j= \{1;2 \}, $

$$
\cov(\xi(i,y_1), \xi(j,y_2)) = {\bf E}\xi(i,y_1)\cdot \xi(j,y_2) = R_{i,j}(y_1,y_2).
$$
 Suppose that the matrix $ R(y_1,y_2) $  for all the values $ Y_1,y_2 \in Y $
 is  less then the following constant matrix
 $ E^{(\rho)} = \{E^{(\rho)}_{i,j}\}, i,j= \{1;2 \}: \ R(y_1,y_2) << E^{(\rho)} $  with entries

 $$
 E^{(\rho)}_{1,1} =  E^{(\rho)}_{2,2} = 1, \  E^{(\rho)}_{1,2} = E^{(\rho)}_{2,1}= \rho,
 $$
where $ \rho = \const \in (-1,1). $ \par
 Recall that the inequality of a view $ A << B $ between two square symmetrical matrices $ A $
and $ B $  with equal size $ l \times l $ is understood as usually: $ A << B $ iff

$$
 \forall x \in R^l \ \Rightarrow \ (Ax,x) \le (Bx,x).
$$
 Assume also the distance $ d $ satisfies the condition (8.5). It follows from the inequality (8.12)
that in the considered case

$$
U(v) \le C(\kappa) \ (v_1^2 - 2 \rho v_1 v_2 + v_2^2 )^{\kappa} \
 \exp \left( - 0.5(1-\rho^2)^{-1}( v_1^2 - 2 \rho v_1 v_2 + v_2^2)  \right), \eqno(8.13)
$$
when $ \min_i v_i \ge 1. $ \par

 Analogous result may be formulated in the multidimensional case. \par

\vspace{3mm}
{\bf Remark 8.1.} More fine result may be obtained by means of the so-called generic chaining method,
see \cite{Ledoux1}, \cite{Talagrand1}, \cite{Talagrand2}, \cite{Talagrand3}, \cite{Talagrand4}.\par

\vspace{3mm}

\section{Concluding remarks}

\vspace{3mm}
{\bf Embedding theorems}. \par
  For the considered in this article m.d.r.i. spaces may be obtained embedding theorems alike
in the one-dimensional case, see \cite{Ostrovsky2}, chapter 1, section 1.17. \\
{\bf Arbitrary measure}. \par
 It may be considered also the case when the measure $ \mu $ is unbounded, but sigma-finite.
The one-dimensional case oh these spaces are investigated, e.g. in  \cite{Ostrovsky3},
\cite{Ostrovsky7}. \par

\end{document}